\documentclass{kjmath}

\usepackage{amssymb}
\usepackage{eucal}
\usepackage{graphicx}
\usepackage{color}
\usepackage{url}

\begin{document}

\kjmvolume{%
@VOL@
}
\kjmnumber{%
@NUM@
}
\kjmyear{%
@YEAR@
}
\kjmstpage{%
@STPAGE@
}
\kjmendpage{%
@ENDPAGE@
}


\title[Novel Constants Based on $\cdots$]{Novel Constants Based on the Generalization of von Neumann-Jordan Constant}


\author[Y.\ Wang]{Yuxin Wang$^1$}
\author[Q.\ Liu]{Qi Liu$^1$}
\author[Q.\ Li]{Qian Li$^1$}

\author[Q.\ Ni]{Qichuan Ni$^1$}

\address{$^1$School of Mathematics and Physics,
\newline \indent Anqing Normal University}
\email{y24060028@stu.aqnu.edu.cn
\newline \indent ORCID iD: \url{https://orcid.org/0009-0003-3091-8728}}

\email{080821074@stu.aqnu.edu.cn
	\newline \indent ORCID iD: \url{https://orcid.org/0009-0008-7101-537X}}
	
\email{liuq67@aqnu.edu.cn
	\newline \indent ORCID iD: \url{https://orcid.org/0000-0002-6049-5282}}
		
\email{y23060030@stu.aqnu.edu.cn
	\newline \indent ORCID iD: \url{https://orcid.org/0009-0003-3091-8728}}

\author[Z.\ Yang]{Zhijian Yang$^2$}

\author[M.\ Sarfraz]{Muhammad Sarfraz$^3$}

\author[Y.\ Li]{Yongjin Li$^3$}

\address{$^2$Department of Mathematics,
\newline \indent Maoming No.1 Middle School}
\email{yangzhj55@mail2.sysu.edu.cn
\newline \indent ORCID iD: \url{https://orcid.org/0009-0001-4669-6221}}

\address{$^3$Department of Mathematics,
	\newline \indent Sun Yat-sen University}
\email{sarfraz@mail2.sysu.edu.cn
\newline \indent ORCID iD: \url{https://orcid.org/0000-0003-4978-9287}}

\email{stslyj@mail.sysu.edu.cn
\newline \indent ORCID iD: \url{https://orcid.org/0000-0003-4322-308X}}

\keywords{Banach spaces, Geometric constants, Normal structure. 
\indent DOI  \\
\indent {\it Received}: \\
\indent {\it Accepted}:}

\begin{abstract}
We introduce a new geometric constant based on a generalization of the parallelogram law, and study its properties as well as some relationships with other well-known geometric constants. A sufficient condition for normal structure is presented. Next, we introduce a constant and calculate its value in a specific space. Furthermore, we introduce another new constant and investigate some of its basic properties.
\end{abstract}

\maketitle

\section{Introduction}
\label{sec:introduction}

In recent years, a lot of geometric constants have been defined and studied in the literature, which makes it easier for us to deal with some problems in Banach space, because it can describe the geometric properties of space quantitatively, and these geometric constants have mathematical beauty, and there are countless relationships between different geometric constants. One of the best known is 
the  von Neumann-Jordan constant $C_{\mathrm{NJ}} (X)$ and  the James constant $J(X)$.
For readers interested in this field are advised to see \cite{JA2,MD,FUY,GJ4,GJ5,JM,NK,MHPP,SDGS}  and the  references mentioned therein. It is worth mentioning that geometric constants play a vital role as a tool for solving other problems, such as in the study of  Banach-Stone theorem, Bishop-Phelps-Bollob{\'a}s theorem, Tingley’s problem, rotundity and smoothness of absolute sums of Banach spaces (see \cite{HA}) and renorming problem in Banach spaces (see \cite{IT}). These are important research topics in the area of functional analysis and we recommend readers to read the literature  \cite{DHC,FC1,RT1}.

Among all normed spaces, the Hilbert spaces are generally considered to have the simplest and clearest geometric structure. Many mathematicians have been finding the conditions on normed spaces under which such spaces become inner product spaces. Results of this kind are of importance in functional analysis, for example, in the theory of operator algebras and certain stability problems  \cite{JC1,PD}.
In addition,  since every physical system is associated with a Hilbert space, the notion of  inner product space also plays a crucial role in quantum mechanics; see   \cite{NZ}.
The rich theory of Hilbert spaces has been created by the efforts of many  mathematicians, we refer
to \cite{DA,AC,MMS} for more details.

The first norm characterization of inner product spaces was given by Fréchet \cite{FS} in 1935. He proved that a normed space $(X, \Vert \cdot \Vert)$ is an inner product space if and only if
\[
\|x+y\|^{2}+\|y+z\|^{2}+\|x+z\|^{2}=\|x+y+z\|^{2}+\|x\|^{2}+\|y\|^{2}+\|z\|^{2}
\]
for all $x,y,z\in X$.

The following great result is  known as ”the parallelogram law”.
\begin{theorem}\rm{\cite{PJJ}}\label{Theorem 2}
	Let $(X, \Vert \cdot \Vert)$ be a real normed linear space. Then $\Vert \cdot \Vert$ derives from an inner product if and only if the parallelogram law holds, i.e.,
	\[2\Vert x\Vert^2+2\Vert y\Vert^2=\Vert x+y\Vert^2+\Vert x-y\Vert^2\]
	for all $x,y \in X$.
\end{theorem}
\begin{theorem}\rm{\cite{MMD}}\label{Theorem 3}
	Let $(X, \Vert \cdot \Vert)$ be a real normed linear space. Then $\Vert \cdot \Vert$ derives from an inner product if and only if
	\[\|x+y\|^{2}+\|x-y\|^{2} \sim 4\]
	for all $x,y \in S_X$, where $\sim$ stands either for $\leq$ or $\geq$.
\end{theorem}

The following theorem is rather well-known  and plays a major role in our article.
For the reader’s convenience, we present the proof as follows:
\begin{theorem}\label{MM1}
	A normed space $(X , \Vert \cdot \Vert)$ is an inner product space if and only if
	\[\lambda\Vert x\Vert^2+(1-\lambda)\Vert y\Vert^2=\Vert \lambda x+(1-\lambda)y\Vert^2+\lambda(1-\lambda)\Vert x-y\Vert^2\]
	for any  $\lambda\in (0,1)$ and $x,y\in X$.
\end{theorem}
\begin{proof} If  $X$ is an inner product space. Observe that for  $\lambda\in (0,1)$ and
	$x,y \in X$  we can obtain
	\[
	\begin{aligned} \lambda\|x\|^{2} &+(1-\lambda)\|y\|^{2}-\lambda(1-\lambda)\|x-y\|^{2} \\ &=\lambda\|x\|^{2}+(1-\lambda)\|y\|^{2}-\lambda(1-\lambda)\left(\|x\|^{2}-2\langle x, y\rangle+\|y\|^{2}\right) \\ &=(\lambda-\lambda(1-\lambda))\|x\|^{2}+2 \lambda(1-\lambda)\langle x, y\rangle+((1-\lambda)-\lambda(1-\lambda))\|y\|^{2} \\ &=\lambda^{2}\|x\|^{2}+2 \lambda(1-\lambda)\langle x, y\rangle+(1-\lambda)^{2}\|y\|^{2} \\ &=\|\lambda x+(1-\lambda) y\|^{2}. \end{aligned}
	\]
	For the second part of the proof, we choose $\lambda=\frac{1}{2}$ and inequality
	\[\lambda\Vert x\Vert^2+(1-\lambda)\Vert y\Vert^2=\Vert \lambda x+(1-\lambda)y\Vert^2+\lambda(1-\lambda)\Vert x-y\Vert^2\]
	becomes
	\[\frac{1}{2}\|x\|^2+\frac{1}{2}\|y\|^2=\bigg\|\frac{1}{2}x+\frac{1}{2}y\bigg\|^2+\frac{1}{4}\|x-y\|^2,\]
	which shows that $X$ is an inner product space.
\end{proof}

Combined with the outstanding work of Jordan and von Neumann \cite{PJJ} on the characterization of inner product spaces with the parallelogram law. Clarkson \cite{JAC1937} first proposed the von Neumann-Jordan constant $C_{\rm{NJ}}(X)$ of Banach space. More precisely, the von Neumann-Jordan constant of $X$ is defined by
\[C_{\mathrm{NJ}}(X)=\sup \bigg\{\frac{\|x+y\|^{2}+\|x-y\|^{2}}{2\left(\|x\|^{2}+\|y\|^{2}\right)}: x, y \in X, (x, y) \neq(0,0)\bigg\},\]
and the famous James constant \cite{JRC3} is defined by
\[
J(X)=\sup\{\min\{\|x+y\|,\|x-y\|\},x,y\in S_X\}.
\]

Moreover, the various axioms of these constants are given in \cite{GJ1,MKYT,MKLMYT}:

${\text {(1)}~ \sqrt2\leq J(X)\leq2}$.

${\text {(2)}~ J(X)=\sqrt 2}~ \text{whenever}~X~\text{represents
	Hilbert space};~\text{the converse is not correct}$.

${\text {(3)}~ 1 \leq C_{\mathrm{NJ}}(X)\leq 2}$.

${\text {(4)}~ X~\text {is a Hilbert space
		iff}~ C_{\mathrm{NJ}}(X)=1}$.

${\text {(5)}~ X \text { is uniformly non-square}~ \text{iff} ~
	C_{\mathrm{NJ}}(X)<2}$.

${\text {(6)}~ C_{\mathrm{NJ}}(X)=C_{\mathrm{NJ}}\left(X^{*}\right)}\
$.

The Modified von Neumann–Jordan constant  
\[
C_{\rm{NJ}}^{\prime}(X)=\sup \left\{\frac{\|x+y\|^{2}+\|x-y\|^{2}}{4}: x, y \in S_{X}\right\}
\]
was studied by Takahashi \cite{TY3} and Alonso et al. \cite{JA2}. It easy to see that
$X$ is uniformly-nonsquare if and only if  $C_{\rm{NJ}}^{\prime}(X)<2$.

Fetter Nathansky et al. \cite{FN} shall introduce new geometric constants: for $\alpha\in[0,1]$,
$$
C_{\alpha}(X)=\sup \bigg\{\frac{\|\alpha x+(1-\alpha)y\|^2+\alpha(1-\alpha)\|x-y\|^2}{\alpha\|x\|^2+(1-\alpha)\|y\|^2}:x,y \in X,\alpha\|x\|^2+(1-\alpha)\|y\|^2 \neq 0 \bigg\},
$$ 
and for $\alpha\in(0,1)$,
$$\begin{aligned}
	D_{\alpha}(X)&=\sup \bigg\{\frac{\|\alpha x+(1-\alpha)y\|^2+\|\alpha x-(1-\alpha)y\|^2+\alpha(1-\alpha)\|x-y\|^2+\alpha(1-\alpha)\|x+y\|^2}{\alpha\|x\|^2+(1-\alpha)\|y\|^2}\\&:x,y \in X, (x,y) \neq (0,0) \bigg\}.
\end{aligned}$$
In addition, some properties of these constants are as follows:

\rm{(1)} $C_{\alpha}(X)=C_{1-\alpha}(X)$ for every $\alpha \in (0,1)$.

\rm{(2)} $C_{\alpha}(X)=1$ for $\alpha \in (0,1)$ if and only if $X$ is a Hilbert space.

\rm{(3)} $C_{\alpha}(X)=C_{\alpha}(X^{*})$ for every $\alpha \in [0,1]$.

\rm{(4)} For $\alpha \in (0,1),D_{1-\alpha}(X)=D_{\alpha}(X) \leq C_{\alpha}(X)$.

\rm{(5)} $1 \leq D_{\alpha}(X) \leq 1+2\sqrt{\alpha(1-\alpha)}$.

\rm{(6)} If $H$ is a Hilbert space,then $D_{\alpha}(H)=1$.

\rm{(7)} $X$ is uniformly nonsquare if and only if $D_{\alpha}(X) \leq 1+2\sqrt{\alpha(1-\alpha)}.$

The Clarkson modulus of convexity \cite{JAC} of a Banach space $X$  is the function
$\delta_X : [0, 2] \rightarrow [0, 1]$ defined by :
\[\delta_X(\epsilon)=\inf \bigg\{1-\frac{\Vert x+y\Vert}{2}:x,y\in S_X, \Vert x-y\Vert\geq \epsilon \bigg\}.\]

The paper is organized as follows:

In section \ref{s1}, we first introduce the constant $C^{\prime}_{\mathrm{NJ}}(\lambda, X)$. The relationship between this new constant and  other well-known constants is  emphasized in terms of nontrivial inequalities. Furthermore, we establish a new 
necessary condition for  Banach spaces have normal structure in the form of $C^{\prime}_{\mathrm{NJ}}(\lambda, X)$. Next, we introduce a constant $C^{\prime\prime}_{\mathrm{NJ}}(\lambda, X)$ and give an example of its value in a particular space. 

In Section \ref{s3}, we introduce a new constant $C^{\prime\prime\prime}_{\mathrm{NJ}}(\lambda, X)$ , this meaning also relates to the  constant $C^{\prime}_{\mathrm{NJ}}(\lambda, X)$ in Section \ref{s1}. Furthermore, some equivalent forms of $C^{\prime\prime\prime}_{\mathrm{NJ}}(\lambda, X)$ and 
its connection  between inner product spaces  are investigated  .

\section{The constant $C^{\prime}_{\mathrm{NJ}}(\lambda, X)$ in the unit sphere}\label{s1}
\label{subsec:CrosRef}

Next we consider the case where $x,y\in S_X$, we shall introduce the constant: for $\lambda\in [0,1]$
\[\begin{aligned}C^{\prime}_{\mathrm{NJ}}(\lambda, X)=\sup\{\Vert \lambda x+(1-\lambda)y\Vert^2+\lambda(1-\lambda)\Vert x-y\Vert^2:x,y\in S_X\}.
\end{aligned}\]

\begin{proposition} \label{p2.1}
Let $X$ be a Banach space. Then
	\[1\leq C^{\prime}_{\mathrm{NJ}}(\lambda, X)\leq4\lambda-4\lambda^2+1.\]
\end{proposition}
\begin{proof}
	Let $x=y\in S_X$, we have
	$$\begin{aligned}
		&\|\lambda x+(1-\lambda)y\|^2+\lambda(1-\lambda)\|x-y\|^2\\ 
		&=\|\lambda x+(1-\lambda)x\|^2+\lambda(1-\lambda)\|x-x\|^2\\
		&\geq \|x\|^2 =1.
	\end{aligned}$$

	The latter assertion can be derived from the following estimate
	$$
	\begin{aligned}
		&\|\lambda x+(1-\lambda)y\|^2+\lambda(1-\lambda)\|x-y\|^2\\ &\leq
		(\lambda\|x\|+(1-\lambda)\|y\|)^2+\lambda(1-\lambda)(\|x\|+\|y\|)^2\\
		&\leq \lambda\|x\|^2+(1-\lambda)\|y\|^2+4\lambda(1-\lambda)\|x\|\|y\|\\
		&\leq 1+4\lambda(1-\lambda)=4\lambda-4\lambda^2+1.
	\end{aligned}$$
\end{proof}

\begin{lemma}\rm{(\cite{AD})}.\label{L1}
{\it Let $X$ be a Banach space.  $\exists \alpha \in [0,\frac{1}{2}]$ so that $\forall u,v \in S_{X}$,
 $\exists \lambda \in [\alpha,1-\alpha]$ with
$$
	\|\lambda u+(1-\lambda)v\Vert^2+\lambda(1-\lambda)\Vert u-v\Vert^2 \sim 1,
$$where $\sim$ stands either for $\leq$ or $\geq$.
Then $X$ is an inner product space.}
\end{lemma}

\begin{proposition}\label{p2.2}
	Let $X$ be a Banach space. 
	Then $C^{\prime}_{\mathrm{NJ}}(\lambda,X)=1$ if and only if $X$ is an inner product space.
\end{proposition}
\begin{proof}
If  $X$ is an inner product space,  using Theorem \ref{MM1} we have
\[\Vert \lambda x+(1-\lambda)y\Vert^2+\lambda(1-\lambda)\Vert x-y\Vert^2=\lambda\Vert x\Vert^2+(1-\lambda)\Vert y\Vert^2,\]
which implies that 
\[\Vert \lambda x+(1-\lambda)y\Vert^2+\lambda(1-\lambda)\Vert x-y\Vert^2\leq \lambda\Vert x\Vert^2+(1-\lambda)\Vert y\Vert^2,\]
for  $\lambda\in (0,1)$ and  $x,y \in S_{X}$.
We have that $C^{\prime}_{\mathrm{NJ}}(\lambda,X)=1.$

For the second part of the proof, if  $C^{\prime}_{\mathrm{NJ}}(\lambda,X)=1$, we have 
$$
\|\lambda u+(1-\lambda)v\Vert^2+\lambda(1-\lambda)\Vert u-v\Vert^2\leq 1.
$$
By the Lemma \ref{L1}, we see that $X$ is an inner product space.
\end{proof}

\begin{proposition}\label{Proposition2}
	For a Banach space $X$, the following assertions are equivalent:
	
	\rm{(i)} $C^{\prime}_{\rm{NJ}}(X)=2$.
	
	\rm{(ii)} $C^{\prime}_{\mathrm{NJ}}(\lambda, X)=4\lambda-4\lambda^2+1$ for all $\lambda\in (0,1)$.
	
	\rm{(iii)} $C^{\prime}_{\mathrm{NJ}}(\lambda, X)=4\lambda_0-4\lambda_0^2+1$ for some $\lambda_0\in (0,1)$.
\end{proposition}
\begin{proof} (i)$\Rightarrow$(ii). Since $C^{\prime}_{\rm{NJ}}(X)=2$, we deduce that thers exists $x_n,y_n\in S_X$ such that
	\[\Vert x_n+y_n\Vert\rightarrow2,~\Vert x_n-y_n\Vert\rightarrow2~(n\rightarrow \infty).\]
	This means that there exists $x_n, y_n \in S_X$ such that
	\[
	\begin{aligned}\|\lambda x_n+(1-\lambda) y_n\| &=\|\lambda(x_n+y_n)-(2\lambda-1) y_n\|\\& \geq \lambda\|x_n+y_n\|-(2\lambda-1)\|y_n\| \\ &=2\lambda-(2\lambda-1)=1.\end{aligned}
	\]

	So we can deduce that there exists $x_n, y_n \in  S_X$ such that
	\[\|\lambda x_n+(1-\lambda) y_n\|\rightarrow1,~ \Vert x_n-y_n\Vert\rightarrow2~
	(n\rightarrow\infty),\]
	which implies that $C^{\prime}_{\mathrm{NJ}}(\lambda, X)= 4\lambda-4\lambda^2+1$.
	
	(ii)$\Rightarrow$(iii). Obvious.
	
	(iii)$\Rightarrow$(i). If $C^{\prime}_{\rm{NJ}}(X)<2$, then there exists $\delta>0$, such
	that for any $x, y \in S_X$ either $\Vert \frac{x+y}{2}\Vert\leq 1-\delta$ or $\Vert \frac{x-y}{2}\Vert\leq 1-\delta$. Without loss of generality, we can assume $\Vert \frac{x-y}{2}\Vert\leq 1-\delta$. 
	Then for some $\lambda_0,\mu_0$ we have
	\[
	\Vert \lambda_0 x+(1-\lambda_0) y\Vert^{2}+\lambda_0
	(1-\lambda_0)\Vert  x- y\Vert^{2}  \leq1+\lambda_0(1-\lambda_0)[2(1-\delta)]^2,
	\]
	which implies that $C^{\prime}_{\mathrm{NJ}}(\lambda, X)<4\lambda_0-4\lambda_0^2+1$.
	This is a contradiction and thus we complete the proof.
\end{proof}

\begin{corollary}\label{Corollary2}
	$X$ is uniformly non-square if and only if  $C^{\prime}_{\mathrm{NJ}}(\lambda, X)<4\lambda-4\lambda^2+1$.
\end{corollary}
\begin{proof}
	It can be directly concluded from Proposition \ref{Proposition2} and 
	the fact $X$ is uniformly-nonsquare if and only if  $C_{\rm{NJ}}^{\prime}(X)<2$.
\end{proof}

\begin{remark}\rm{Corollary} \ref{Corollary2}  can also be interpreted in another way.
	It is clear that
	\[\begin{aligned} &\min \{ \Vert \lambda x+(1-\lambda)y\Vert,\sqrt{\lambda(1-\lambda)}\Vert x-y\Vert\} \\
		&\leq \bigg(\frac {1}{2} (\Vert \lambda x+(1-\lambda)y\Vert^2+\lambda(1-\lambda)\Vert x-y\Vert^2)\bigg)^{\frac {1}{2}}
		\\& \leq \bigg(\frac{1}{2}C^{\prime}_{\mathrm{NJ}}(\lambda, X) (\lambda\Vert x\Vert^2+(1-\lambda)\Vert y\Vert^2)\bigg)^{\frac{1}{2}}\end{aligned}\]
	and hence
	\[ \min \{ \Vert \lambda x+(1-\lambda)y\Vert,\sqrt{\lambda(1-\lambda)}\Vert x-y\Vert\}\leq\sqrt{\frac{1}{2}C^{\prime}_{\mathrm{NJ}}(\lambda, X)}.
	\]
	for any $x, y\in S_X$.
	Then, by  setting  $\lambda=\frac{1}{2}$ leads to the estimates
	\[\begin{aligned}\min \{ \| x +  y\|, \|x -  y\|\}&\leq 2\sqrt{\frac{1}{2}C^{\prime}_{\mathrm{NJ}}(\frac{1}{2},X)}\\&=2(1-\epsilon),\end{aligned}\]
	where \[\epsilon=1-\sqrt{\frac{1}{2}C^{\prime}_{\mathrm{NJ}}(\frac{1}{2},X)}.\]
	Using  $C^{\prime}_{\mathrm{NJ}}(\lambda, X)<4\lambda-4\lambda^2+1$  it follows that  $0<\epsilon<1,$
	which shows that
	\[\min \{ \| x +  y\|, \|x -  y\|\}\leq 2(1-\epsilon),\]
	as stated.
\end{remark}

The following proposition characterizes inner product spaces using convex functions. For more information on the relationship between convex functions and inner product spaces, please refer to \cite{NP}.

\begin{proposition} \label{Proposition3.6}
	Let $X$ be a Banach space.  Then
	
	\rm{(i)}  $C^{\prime}_{\mathrm{NJ}}(\lambda, X)$ is a convex function if and only if $X$ is an inner product space.
	
	\rm{(ii)} $C^{\prime}_{\mathrm{NJ}}(\lambda, X)$ is continuous on $(0,1)$.
\end{proposition}
\begin{proof}
	\rm{(i)} If  $X$ is an inner product space,  using Theorem \ref{MM1} we have
	\[\Vert \lambda x+(1-\lambda)y\Vert^2+\lambda(1-\lambda)\Vert x-y\Vert^2=\lambda\Vert x\Vert^2+(1-\lambda)\Vert y\Vert^2=1\]
	for $\lambda\in (0,1)$ and  $x,y \in S_X$, which implies that 
	$C^{\prime}_{\mathrm{NJ}}(\lambda, X)\equiv1$, as desired.

	For the second part of the proof, firstly we have $C^{\prime}_{\mathrm{NJ}}(\lambda, X)=C^{\prime}_{\mathrm{NJ}}(1-\lambda, X)$,which implies that 
	\[C^{\prime}_{\mathrm{NJ}}(0,X)=C^{\prime}_{\mathrm{NJ}}(1,X)=1.\]
	
	Therefore, from   the fact $C^{\prime}_{\mathrm{NJ}}(\lambda, X)$ is 
	convex function, we obtain
	\[\begin{aligned}1&\leq  C^{\prime}_{\mathrm{NJ}}(\lambda, X)\\&=
		C^{\prime}_{\mathrm{NJ}}((1-\lambda)\cdot0+\lambda\cdot1,X)
		\\&\leq (1-\lambda)C^{\prime}_{\mathrm{NJ}}(0,X)+\lambda C^{\prime}_{\mathrm{NJ}}(1,X)\\&=1
	\end{aligned}\]
	for $\lambda\in (0,1)$, and hence $C^{\prime}_{\mathrm{NJ}}(\lambda,X)=1$.
	This implies that $X$ is an inner product space.

	\rm{(ii)} We distinguish two cases.
	
	{\bf Case 1.} $\forall \lambda \in (\frac{1}{2},1)$.
	
	$\forall\epsilon>0$, $\exists\delta_0=\min\{\frac{\epsilon
	}{8},\frac{1-\lambda}{2}\}>0$.  When $0<\delta<\delta_0$,  we conclude that the following estimate holds.
	\[\begin{aligned}
		&\|(\lambda+\delta) x+(1-(\lambda+\delta)) y\|^{2}+(\lambda+\delta)(1-(\lambda+\delta))\|x-y\|^{2}\\
		&-\|\lambda x+(1-\lambda) y\|^{2}-\lambda(1-\lambda)\|x-y\|^{2}
		\\&=\|(\lambda+\delta) x+(1-(\lambda+\delta)) y\|^{2}-\|\lambda x+(1-\lambda) y\|^{2}+\delta(1-2 \lambda-\delta)\|x-y\|^{2} 
		\\&\leq(\|(\lambda+\delta) x+(1-(\lambda+\delta)) y\|-\|\lambda x+(1-\lambda) y\|)(\|(\lambda+\delta) x+(1-(\lambda+\delta)) y\|\\
		&+\|\lambda x+(1-\lambda) y\|)  
		\\&\leq\|(\lambda+\delta) x+(1-(\lambda+\delta)) y-(\lambda x+(1-\lambda) y)\|((\lambda+\delta)\|x\|\\
		&+(1-(\lambda+\delta))\|y\|+\lambda\|x\|+(1-\lambda)\|y\|)
		\\&=2\delta\Vert x-y\Vert\leq4\delta<\epsilon.
	\end{aligned}\]
	Then we obtain
	\[
	\begin{aligned}&\bigg|C^{\prime}_{\mathrm{NJ}}(\lambda+\delta, X)-C^{\prime}_{\mathrm{NJ}}(\lambda, X)\bigg| \\& =\bigg| \sup \left\{\|(\lambda+\delta) x+(1-(\lambda+\delta)) y\|^{2}+(\lambda+\delta)(1-(\lambda+\delta))\|x-y\|^{2}: x, y\in S_X\right\} \\& \quad-\sup \left\{\|\lambda x+(1-\lambda) y\|^{2}+\lambda(1-\lambda)\|x-y\|^{2}: x, y \in S_{X}\right\} \bigg| \\&\leq
		\bigg|\sup \{\|(\lambda+\delta) x+(1-(\lambda+\delta)) y\|^{2}+(\lambda+\delta)(1-(\lambda+\delta))\|x-y\|^{2} \\& \quad-\|\lambda x+(1-\lambda) y\|^{2}-\lambda(1-\lambda)\|x-y\|^{2}:x,y\in S_X\}\bigg|
		<\epsilon,\end{aligned}
	\]
	which implies that 
	\[\bigg|C^{\prime}_{\mathrm{NJ}}(\lambda+\delta, X)-C^{\prime}_{\mathrm{NJ}}(\lambda, X)\bigg|<\epsilon.\]
	
	{\bf Case 2.} $\forall \lambda \in (0,\frac{1}{2}]$.
	
	Let $\bar{\lambda}=1-\lambda\in [\frac{1}{2},1)$. It is clear that
	$C^{\prime}_{\mathrm{NJ}}(\bar{\lambda},X)$ is continuous on  $[\frac{1}{2},1)$. This means that for
	every $\epsilon>0$ there exists a $\delta_0>0$ such that
	\[\bigg|C^{\prime}_{\mathrm{NJ}}(\bar{\lambda_2}, X)-C^{\prime}_{\mathrm{NJ}}(\bar{\lambda_1}, X)\bigg|<\epsilon~
	\operatorname{if}~
	|\bar{\lambda_2}-\bar{\lambda_1}|<\delta_0, \bar{\lambda_1},\bar{\lambda_2}\in \big[\frac{1}{2},1).\]
	Note also that
	\[C^{\prime}_{\mathrm{NJ}}(\lambda, X)=C^{\prime}_{\mathrm{NJ}}(1-\lambda, X)=C^{\prime}_{\mathrm{NJ}}(\bar{\lambda}, X),\]
	we get 
	\[\bigg|C^{\prime}_{\mathrm{NJ}}(\lambda_2, X)-C^{\prime}_{\mathrm{NJ}}(\lambda_1, X)\bigg|=
	\bigg|C^{\prime}_{\mathrm{NJ}}(\bar{\lambda_2}, X)-C^{\prime}_{\mathrm{NJ}}(\bar{\lambda_1}, X)\bigg|<\epsilon.\]
\end{proof}

\begin{example}
	Consider $X=l_3$, then $C^{\prime}_{\mathrm{NJ}}(\frac{1}{2},X)$  be  not a convex function.
	Indeed, take \[x_0=(1,0,\ldots), y_0=\bigg(\frac{1}{\sqrt[3]{2}},\frac{1}{\sqrt[3]{2}},0,\ldots\bigg).\]
	It is clear that $\Vert x_0\Vert_3=\Vert y_0\Vert_3=1$, and
	\[C^{\prime}_{\mathrm{NJ}}(\frac{1}{2},X)\geq \bigg\Vert \frac{x_0+y_0}{2}\bigg\Vert^2+
	\frac{1}{4}\Vert x_0-y_0\Vert^2.\]
	By elementary calculus, we can easily get 
	\[\bigg\Vert \frac{x_0+y_0}{2}\bigg\Vert^2+
	\frac{1}{4}\Vert x_0-y_0\Vert^2>1\]
	and hence 
	\[C^{\prime}_{\mathrm{NJ}}(\frac{1}{2}, X)>1=\frac{1}{2}C^{\prime}_{\mathrm{NJ}}(0, X)+\frac{1}{2}C^{\prime}_{\mathrm{NJ}}(1, X).\]
	
\end{example}

\begin{theorem}
	Let $X$ be a Banach space, $\epsilon\in [0,2]$.
	Then 
	\[ \begin{aligned}\frac{1}{4}\epsilon^2+|2\lambda-1|(1-\delta_X(\epsilon))\epsilon
		+(1-\delta_X(\epsilon))^2&\leq C^{\prime}_{\mathrm{NJ}}(\lambda,X)\\&\leq (2\lambda+|1-2\lambda|-2\lambda\delta_X(\epsilon))^2+\lambda(1-\lambda)\epsilon^2.\end{aligned}\]
\end{theorem}
\begin{proof} 
	Let $\delta_X(\epsilon)=\alpha$. We can deduce that $\exists~ x_n,y_n\in S_X, \Vert x_n-y_n\Vert=\epsilon$,
	such that 
	$
	\lim_{n \to \infty} \Vert x_n+y_n\Vert=2(1-\alpha)
	$.
	Furthermore, 
	\[\begin{aligned}\Vert \lambda x_n+(1-\lambda)y_n\Vert&=
		\bigg\Vert \frac{1}{2}(x_n+y_n)+(\lambda-\frac{1}{2})(x_n-y_n)\bigg\Vert
		\\&\geq \bigg| \frac{1}{2}\Vert x_n+y_n\Vert-|\lambda-\frac{1}{2}|\Vert x_n-y_n\Vert\bigg|
		\\&=\frac{1}{2}\bigg|\Vert x_n+y_n\Vert-|2\lambda-1|\epsilon \bigg|,
	\end{aligned}\]
	which implies that
	\[\begin{aligned}C^{\prime}_{\mathrm{NJ}}(\lambda,X)&\geq \frac{1}{4}\bigg| \Vert 
		x_n+y_n\Vert-|2\lambda-1|\epsilon\bigg|^2+\lambda(1-\lambda)\epsilon^2.\end{aligned}\]
	Let $n\rightarrow \infty$, as desired.

	On the other hand, we can deduce that the following estimate holds.
	\[\begin{aligned}\Vert \lambda x+(1-\lambda) y\Vert&
		\leq \lambda \Vert x+y\Vert+|1-2\lambda|\Vert y\Vert\\&
		\leq \lambda(2-2\delta(\epsilon))+|1-2\lambda|\\& =2\lambda+|1-2\lambda|-2\lambda\delta(\epsilon).
	\end{aligned}\]
	Then we obtain 
	\[\Vert \lambda x+(1-\lambda)y\Vert^2+\lambda(1-\lambda)\Vert x-y\Vert^2\leq 
	(2\lambda+|1-2\lambda|-2\lambda\delta_X(\epsilon))^2+\lambda(1-\lambda)\epsilon^2.\]
	This completes the proof.
\end{proof}

\begin{example}
	Consider $X$ be $c_{0}=\{\left(x_{i}\right):\displaystyle\lim _{i \rightarrow \infty}|x_{i}|=0\}$ be equipped with the norm defined by
	\[
	\|x\|=\sup _{1 \leq i<\infty}\left|x_{i}\right|+\left(\sum_{i=1}^{\infty} \frac{\left|x_{i}\right|^{2}}{4^{i}}\right)^{1 / 2}.
	\]
	It is known that $\delta_X(2)=1$(see \cite{UA}). 
	Then we have
	\[\begin{aligned}  C^{\prime}_{\mathrm{NJ}}(\lambda,X)
		&\leq(2\lambda+|1-2\lambda|-2\lambda\delta_X(\epsilon))^2+\lambda(1-\lambda)\epsilon^2
		\\&=1.\end{aligned}\]
	Thus $C^{\prime}_{\mathrm{NJ}}(\lambda,X)=1$ directly from the estimate of $C^{\prime}_{\mathrm{NJ}}(\lambda,X)$.
\end{example}

Next, we will see that the  normal structure and
the constant $C^{\prime}_{\mathrm{NJ}}$ have  a close relationship. 
\begin{definition}\rm{\cite{MSB}}
	Let $K$ be a non-singleton subset of a Banach space $X$, if $K$ is
	closed, bounded as well as convex then $X$ holds the normal
	structure, whenever $r(K)<\operatorname{diam}(K)$ for every $K$,
	where $r(K)$ and $\operatorname{diam}(K)$ are respectively
	symbolized for diameter as well as for Chebyshev radius, and
	consequently defined mathematically as is
	\[\operatorname{diam}(K):=\sup \{\Vert x-y\Vert:x,y\in K\}\]
	and
	\[r(K):=\inf\{\sup\{\Vert x-y\Vert: y\in K\}
	:x\in K\}.\]
\end{definition}
Normal structure is an important concept in fixed point theory \cite{SP}. A Banach space $X$ is said to have weak normal structure if each weakly compact convex set $K$ of $X$ that contains more than one point has normal structure. Even more to the point, for Banach space $X$ which is reflexive, the weak normal structure and  normal structure coincide. Furthermore, every reflexive Banach space with normal structure has the fixed point property.

We begin by starting a lemma which will be our main tool.
\begin{lemma}\rm{(\cite{GJ1})}.\label{1} {\it Let $X$ be a Banach space without weak normal structure, then for any $0 < \delta < 1$, there exist $x_1, x_2, x_3$ in $S_X$ satisfying
		
		{\rm (i)}  $x_2-x_3 =ax_1$ with $|a-1|<\delta$;
		
		{\rm (ii)}  $|\Vert x_1 - x_2\Vert - 1|, |\Vert x_3 - (-x_1)\Vert - 1| < \delta$; and
		
		{\rm (iii)} $\Vert \frac{x_1+x_2}{2}\Vert,\Vert \frac{x_3+(-x_1)}{2}\Vert>1-\delta$.}
\end{lemma}

It is possible to understand the geometric sense of this Lemma as follows: if $X$ does not have a weak normal structure, then an inscribed hexagon exists in $S_X$ with an arbitrarily closed length of each side to $1$, and at least four sides with an arbitrarily small distance to $S_X$.

\begin{theorem}
	A Banach space $X$ with $C^{\prime}_{\mathrm{NJ}}(\lambda,X)<-\lambda^2+\lambda+1$ for some $\lambda\in [\frac{1}{2},1)$ has  normal structure.
\end{theorem}
\begin{proof} Notice that, by  Corollary \ref{Corollary2},
	$C^{\prime}_{\mathrm{NJ}}(\lambda,X)<-\lambda^2+\lambda+1$  
	implies that $X$ is uniformly non-square, and hence reflexive \cite{JRC3}. This means that normal structure and weak normal structure coincide.

	Suppose $X$ does not have weak normal structure.   For
	each $\delta > 0$, let $x_1,x_2$ and $x_3$ in $S_X$ satisfying the conditions in
	Lemma \ref{1}. 
	
	Then
	\[\begin{aligned}
		\bigg\| x_1+\frac{1-\lambda}{\lambda}x_2\bigg\|&=\bigg\| (x_1+x_2)-(1-\frac{1-\lambda}{\lambda})x_2\bigg\|\\&\geq
		\Vert x_1+x_2\Vert-\bigg\|(1-\frac{1-\lambda}{\lambda})x_2\bigg\|\\&\geq
		2-2\delta-\bigg(1-\frac{1-\lambda}{\lambda}\bigg)
	\end{aligned}\]
	
	and
	\[\lambda(1-\lambda)\Vert x_1-x_2\Vert^2\geq\lambda(1-\lambda)(1-\delta)^2.
	\]
	
	It is easy to see that $\delta$ can be arbitrarily small,  so we obtain
	\[\begin{aligned}C^{\prime}_{\mathrm{NJ}}(\lambda,X)&\geq
		\lambda^2\Vert  x+\frac{(1-\lambda)}{\lambda}y\Vert^2+\lambda(1-\lambda)\Vert x-y\Vert^2\\&\geq
		-\lambda^2+\lambda+1,\end{aligned}\]
	which is a contradiction. This completes the proof. 
\end{proof}

Next, we consider the special case where $x$ and $y$ on the unit sphere satisfy isosceles orthogonality, we  shall introduce the constant: for $\lambda\in [0,1]$
\[\begin{aligned}C^{\prime\prime}_{\mathrm{NJ}}(\lambda,X)=\sup\{\Vert \lambda x+(1-\lambda)y\Vert^2+\lambda(1-\lambda)\Vert x-y\Vert^2:x,y\in S_X,x\perp_{I}y\}.
\end{aligned}\]

\begin{remark}
	Clearly, $C^{\prime\prime}_{\mathrm{NJ}}(\lambda,X)\leq C^{\prime}_{\mathrm{NJ}}(\lambda,X)$ is always valid for any Banach spaces.
\end{remark}
 
 The next proposition will present the upper and lower bound estimates of $C^{\prime\prime}_{\mathrm{NJ}}(\lambda,X)$.
\begin{proposition} Suppose that $X$ is a Banach space. Then
	\[1\leq C^{\prime\prime}_{\mathrm{NJ}}(\lambda,X)\leq4\lambda-4\lambda^2+1.\]
\end{proposition}
\begin{proof}
	By using the same method as Proposition \ref{p2.1}, that we can easily get the result, so we omit the proof.
\end{proof}

\begin{example}
	Let $\lambda\in[0,\frac{1}{3}]$ or  $\lambda\in[\frac{1}{3},1]$, $X=\mathbb{R}^2$ with $\ell_{\infty}-\ell_1$ norm defined by   
	\[
	\|x\|=\left\{\begin{array}{l}
		\|x\|_1, x_1 x_2 \leq 0, \\
		\|x\|_{\infty}, x_1 x_2 \geq 0.
	\end{array}\right.
	\]
	Then  \[
	C^{\prime\prime}_{\mathrm{NJ}}(\lambda,X) =1.2.
	\]
\end{example}

First, consider  $\lambda\in[0,\frac{1}{3}]$.

If $x= (y_1,1+y_1 ),y= (y_2,1 +y_2)$, where $-1\leq y_1\leq y_2\leq0$; $x=(y_1,y_1-1),y=(y_2,y_2-1)$, where $0\leq y_1\leq y_2\leq 1$. For both cases, it is determined by $x \perp_{I} y$, there  $|y_1-y_2|=2$, which is contradictory. To calculate the value of the constant $C^{\prime\prime}_{\mathrm{NJ}}(\lambda,X)$, there are only two cases to consider.

{\bf 
	Case 1:} Assuming that  $x=\left(x_1, 1\right), y=\left(1, y_2\right)$,  $0 \leq x_1 \leq y_2 \leq 1$. Since  $x \perp_I y$, we have 
\[
1+y_2=\left(1-x_1\right)+\left(1-y_2\right) \text {, }
\]
hence $x_1+2 y_2=1$,   $y_2 \in\left[\frac{1}{3}, \frac{1}{2}\right]$.

Then, $\|\lambda x+(1-\lambda)y\|=1-2\lambda y_2, \|x-y\|=1+y_2.$

By the definition of the constant $C^{\prime\prime}_{\mathrm{NJ}}(\lambda,X)$, we have
\[C^{\prime\prime}_{\mathrm{NJ}}(\lambda,X)=\max_{\frac{1}{3}\leq y_2 \leq \frac{1}{2}}(1-2\lambda y_2)^2+\lambda(1-\lambda)(1+y_2)^2.\]

By simple calculation, we get $C^{\prime\prime}_{\mathrm{NJ}}(\lambda,X)=1.04$.

{\bf 
	Case 2:} Assuming that  $x=\left(x_1, 1\right), y=\left(y_1, 1+y_1\right)$ satisfies $-1 \leq y_1 \leq 0 \leq x_1 \leq 1$. Since $x\perp_{I} y$, we have $\|
(x_1+y_1,2+y_1)\|=\|(x_1-y_1,-y_1)\|$. Now, we discuss cases (i) and (ii).

(i) If $-x_1 \leq y_1$, then $2+y_1=x_1-y_1$ is true, hence $x_1-2y_1=2$,  $y_1 \in\left[-\frac{2}{3},-\frac{1}{2}\right].$ 

Then $\|\lambda x+(1-\lambda)y\|=1-2\lambda y_1-2\lambda, \|x-y\|=2+y_1.$

By the definition of the constant $C^{\prime\prime}_{\mathrm{NJ}}(\lambda,X)$, we have
\[C^{\prime\prime}_{\mathrm{NJ}}(\lambda,X)=\max_{-\frac{2}{3}\leq y_1 \leq-\frac{1}{2}}\left\{(1-2\lambda y_1-2\lambda)^2+\lambda(1-\lambda)(2+y_1)^2\right\}.\]

By simple calculation, we get $C^{\prime\prime}_{\mathrm{NJ}}(\lambda,X)=1.04$.

(ii) If $y_1 \leq-x_1$, then $(-x_1-y_1)+(2+y_1)=x_1-y_1$ is true, hence $2x_1=2+y_1$,  $y_1 \in\left[-1,-\frac{2}{3}\right].$ 

Then $\|\lambda x+(1-\lambda)y\|=1-\lambda -\frac{\lambda}{2}y_1, \|x-y\|=1-\frac{y_1}{2}.$

By the definition of the constant $C^{\prime\prime}_{\mathrm{NJ}}(\lambda,X)$, we have
\[C^{\prime\prime}_{\mathrm{NJ}}(\lambda,X)=\max_{-1\leq y_1 \leq-\frac{2}{3}}\left\{(1-\lambda -\frac{\lambda}{2}y_1)^2+\lambda(1-\lambda)\bigg(1-\frac{y_1}{2}\bigg)^2\right\}.\]

By simple calculation, we get $C^{\prime\prime}_{\mathrm{NJ}}(\lambda,X)=1.2$.

Combining {\bf Case 1} and {\bf Case 2}, we have $C^{\prime\prime}_{\mathrm{NJ}}(\lambda,X)=1.2$.

As for the case of  $\lambda\in[\frac{1}{3},1]$, we can do a similar discussion, omit the process, and still get $C^{\prime\prime}_{\mathrm{NJ}}(\lambda,X)=1.2$.

\section{ Properties of the constant $C^{\prime\prime\prime}_{\mathrm{NJ}}(\lambda,X)$}\label{s3}
\label{subsec:SecUnit}
From now on, we will consider only Banach spaces of dimension at least $2$. We begin by introducing the following key definition: for $\lambda\in [0,1]$
$$\begin{aligned}
C^{\prime\prime\prime}_{\mathrm{NJ}}(\lambda,X)&=\sup \bigg\{\frac{1}{2}(\|\lambda x+(1-\lambda)y\|^2+\|\lambda x-(1-\lambda)y\|^2+\lambda(1-\lambda)\|x-y\|^2\\&+\lambda(1-\lambda)\|x+y\|^2): x,y \in S_{X}\bigg\}.
\end{aligned}$$

\begin{proposition}\label{Proposition2.1} Let $X$ be a Banach space. Then
	$$
	1 \leq C^{\prime\prime\prime}_{\mathrm{NJ}}(\lambda,X) \leq -4\lambda^2+4\lambda+1.
	$$
\end{proposition}

\begin{proof} 
	Let $x=y\in S_X$, we have
	$$\begin{aligned}
		&\frac{1}{2}(\|\lambda x+(1-\lambda)y\|^2+\|\lambda x-(1-\lambda)y\|^2+\lambda(1-\lambda)\|x-y\|^2+\lambda(1-\lambda)\|x+y\|^2)\\ &\ge \frac{1}{2}(\|x\|^2+(2\lambda-1)^2\|x\|^2+4\lambda(1-\lambda)\|x\|^2) =1,
	\end{aligned}$$which means that  $C^{\prime\prime\prime}_{\mathrm{NJ}}(\lambda,X)\geq 1$.

	The latter assertion can be derived from the following estimate
	$$
	\begin{aligned}
		&\frac{1}{2}(\|\lambda x+(1-\lambda)y\|^2+\|\lambda x-(1-\lambda)y\|^2+\lambda(1-\lambda)\|x-y\|^2+\lambda(1-\lambda)\|x+y\|^2)\\ &\le
		\frac{1}{2}((\lambda\|x\|+(1-\lambda)\|y\|)^2+(\lambda\|x\|+(1-\lambda)\|y\|)^2+\lambda(1-\lambda)(\|x\|+\|y\|)^2\\&+\lambda(1-\lambda)(\|x\|+\|y\|)^2)\\
		&=\lambda\|x\|^2+(1-\lambda)\|y\|^2+4\lambda(1-\lambda)\|x\|\|y\|\\
		&\leq 1+4\lambda(1-\lambda)= -4\lambda^2+4\lambda+1.
	\end{aligned}$$
	This implies that $C^{\prime\prime\prime}_{\mathrm{NJ}}(\lambda,X) \leq -4\lambda^2+4\lambda+1,$ as desired.
\end{proof}

\begin{proposition}\label{Proposition3} Let $X$ be a Banach space. Then
	\[ C^{\prime\prime\prime}_{\mathrm{NJ}}(\lambda,X)=C^{\prime\prime\prime}_{\mathrm{NJ}}(1-\lambda,X).\]
\end{proposition}
\begin{proof} By replacing $\lambda$ with $1-\lambda$,we can easily get the result,so we omit the proof. 
\end{proof}

We now give the following Proposition, which is inspired by Theorem \ref{MM1}.

\begin{proposition}\label{MM2}
	Let  $(X , \Vert \cdot \Vert)$ be a normed space. Then $\|\cdot\|$ derives from an
	inner product if and only if
	$$\begin{aligned}
		&\|\lambda x+(1-\lambda)y\|^2+\|\lambda x-(1-\lambda)y\|^2+\lambda(1-\lambda)\|x-y\|^2+\lambda(1-\lambda)\|x+y\|^2\\& \leq 2(\lambda\|x\|^2+(1-\lambda)\|y\|^2).
	\end{aligned}$$
	for any non-negative real numbers $\lambda\in (0,1)$ and  $x,y \in S_{X}$.
\end{proposition}

\begin{proof} If  $X$ is an inner product space,  using Theorem \ref{MM1} we have
	\[\Vert \lambda x+(1-\lambda)y\Vert^2+\lambda(1-\lambda)\Vert x-y\Vert^2=\lambda\Vert x\Vert^2+(1-\lambda)\Vert y\Vert^2,\]
	\[\Vert \lambda x-(1-\lambda)y\Vert^2+\lambda(1-\lambda)\Vert x+y\Vert^2=\lambda\Vert x\Vert^2+(1-\lambda)\Vert y\Vert^2.\]
	for  $\lambda\in (0,1)$ and  $x,y \in S_{X}$.
	
	For the second part of the proof, we choose $\lambda=\frac{1}{2}$ and inequality
	\[\lambda\Vert x\Vert^2+(1-\lambda)\Vert y\Vert^2\leq\Vert \lambda x+(1-\lambda)y\Vert^2+\lambda(1-\lambda)\Vert x-y\Vert^2,\]
	\[\lambda\Vert x\Vert^2+(1-\lambda)\Vert y\Vert^2\leq\Vert \lambda x-(1-\lambda)y\Vert^2+\lambda(1-\lambda)\Vert x+y\Vert^2,\]
	becomes
	\[\frac{1}{2}\|x\|^2+\frac{1}{2}\|y\|^2\leq\bigg\|\frac{1}{2}x+\frac{1}{2}y\bigg\|^2+\frac{1}{4}\|x-y\|^2,\]
	\[\frac{1}{2}\|x\|^2+\frac{1}{2}\|y\|^2\leq\bigg\|\frac{1}{2}x-\frac{1}{2}y\bigg\|^2+\frac{1}{4}\|x+y\|^2.\]
	By Theorem \ref{Theorem 3}, we see that $X$ is an inner product space for any non-negative real numbers $\lambda \in (0,1)$.
\end{proof}

\begin{theorem}Let $X$ be a Banach space. 
	Then $C^{\prime\prime\prime}_{\mathrm{NJ}}(\lambda,X)=1$ if and only if $X$ is an inner product space, for any non-negative real numbers $\lambda \in (0,1)$.
\end{theorem}
\begin{proof} By using the same method as Proposition \ref{p2.2}, that we can easily get the result,
	so we omit the proof.
\end{proof}

\begin{proposition}
	Let $X$ be a Banach space. Then,
	$$\begin{aligned}
		2\min\{\lambda,1-\lambda\}C^{'}_{NJ}(X) &\leq C^{\prime\prime\prime}_{\mathrm{NJ}}(\lambda,X)\\ &\leq 2\max\{\lambda,1-\lambda\}C^{'}_{NJ}(X)+2\sqrt{2}\max\{\lambda,1-\lambda\}\\&|\lambda-(1-\lambda)|\sqrt{C^{'}_{NJ}(X)}+|\lambda-(1-\lambda)|^2.
	\end{aligned}$$
\end{proposition}
\begin{proof} 
	According to the following elementary identity:
	$$
	\lambda x+(1-\lambda)y=\frac{\lambda+(1-\lambda)}{2}(x+y)+\frac{\lambda-(1-\lambda)}{2}(x-y),
	$$
	$$
	\lambda x-(1-\lambda)y=\frac{\lambda-(1-\lambda)}{2}(x+y)+\frac{\lambda+(1-\lambda)}{2}(x-y).
	$$
	Then, we can obtain
	$$\begin{aligned}
		C^{\prime\prime\prime}_{\mathrm{NJ}}(\lambda,X) &\geq \frac{1}{2}\bigg(\bigg(\bigg\|\frac{\lambda+(1-\lambda)}{2}(x+y)\bigg\|-\bigg\|\frac{\lambda-(1-\lambda)}{2}(x-y)\bigg\|\bigg)^2+\bigg(\bigg\|\frac{\lambda-(1-\lambda)}{2}(x+y)\bigg\|\\&-\bigg\|\frac{\lambda+(1-\lambda)}{2}(x-y)\bigg\|\bigg)^2
		+\lambda(1-\lambda)\|x-y\|^2+\lambda(1-\lambda)\|x+y\|^2\bigg)\\
		&\geq \frac{1}{2}\bigg(\frac{1}{2}\|x-y\|^2+\frac{1}{2}\|x+y\|^2-(\lambda^2-(1-\lambda)^2\|x-y\|\|x+y\|\bigg),
	\end{aligned}$$
	which implies that
	$$\begin{aligned}
		2C^{\prime\prime\prime}_{\mathrm{NJ}}(\lambda,X)&\geq \frac{1}{2}\|x-y\|^2+\frac{1}{2}\|x+y\|^2-(\lambda+(1-\lambda))|\lambda-(1-\lambda)|\|x-y\|\|x+y\|\\&
		\geq \frac{1}{2}\|x-y\|^2+\frac{1}{2}\|x+y\|^2-\frac{1}{2}(\lambda+(1-\lambda))|\lambda-(1-\lambda)|(\|x-y\|^2+\|x+y\|^2).
	\end{aligned}$$
	This shows that
	$$
	C^{\prime\prime\prime}_{\mathrm{NJ}}(\lambda,X) \geq 2\min\{\lambda,1-\lambda\}C^{'}_{NJ}(X).
	$$
	On the other hand, we can deduce that 
	$$
	\begin{aligned}
		C^{\prime\prime\prime}_{\mathrm{NJ}}(\lambda,X) &\leq \frac{1}{2}((\lambda\|(x+y)\|+|\lambda-(1-\lambda)|\|y\|)^2+((1-\lambda)\|x-y\|+|\lambda-(1-\lambda)|\|x\|)^2\\
		&+\lambda(1-\lambda)\|x-y\|^2+\lambda(1-\lambda)\|x+y\|^2)\\
		&\leq \frac{1}{2}(\max\{\lambda,1-\lambda\}(\|x-y\|^2+\|x+y\|^2)+|\lambda-(1-\lambda)|^2(\|x\|^2+\|y\|^2)\\
		&+2\max\{\lambda,1-\lambda\}|\lambda-(1-\lambda)|\sqrt{\|x-y\|^2+\|x+y\|^2}\sqrt{\|x\|^2+\|y\|^2})\\
		&\leq 2\max\{\lambda,1-\lambda\}C^{'}_{NJ}(X)+2\sqrt{2}\max\{\lambda,1-\lambda\}|\lambda-(1-\lambda)|\sqrt{C^{'}_{NJ}(X)}\\&+|\lambda-(1-\lambda)|^2.
	\end{aligned}
	$$
\end{proof}

\begin{theorem}
	Let $X$ be a Banach space, $\epsilon\in [0,2]$. If 
	\[C^{\prime\prime\prime}_{\mathrm{NJ}}(\lambda,X)<\frac{1}{2}\lambda \epsilon^{2}-(\epsilon-4)\lambda+1,\]
	then $\delta_X(\epsilon)>0$.
\end{theorem}
\begin{proof} Suppose $\delta_X(\epsilon)=0$, then there exist $x_n,y_n\in S_X$ such
	that $\Vert x_n-y_n\Vert=\epsilon$ for all $n\in \mathbb{N}$ and $\displaystyle\lim_{n\rightarrow\infty}
	\Vert x_n+y_n\Vert=2$.
	According to the following elementary inequality:
	$$
	\|\lambda x_{n}+(1-\lambda) y_{n}\| \geq|\|\lambda(x_{n}-y_{n})\|-\|y_{n}\||=|\lambda \epsilon-1|,
	$$
	$$
	\|\lambda x_{n}-(1-\lambda) y_{n}\| \geq|\|\lambda(x_{n}+y_{n})\|+\|y_{n}\||=|2\lambda+1|,
	$$
	for any $x_n,y_n\in S_X$.
	
	We can deduce that
	$$\begin{aligned}
		\frac{1}{2}\lambda \epsilon^{2}-(\epsilon-4)\lambda+1&=\frac{1}{2}(\lambda \epsilon^{2}-(2\epsilon-8)\alpha+2)\\
		&=\frac{1}{2}((\lambda \epsilon-1)^2+(2\lambda+1)^2+\lambda(1-\lambda)\epsilon^2+4\lambda(1-\lambda))\\
		&\leq \lim_{n \rightarrow 0}\inf \frac{1}{2}(\|\lambda x_{n}+(1-\lambda)y_{n}\|^2+\|\lambda x_{n}-(1-\lambda)y_{n}\|^2\\&+\lambda(1-\lambda)\|x_{n}-y_{n}\|^2
		+\lambda(1-\lambda)\|x_{n}+y_{n}\|^2)\\
		&\leq C^{\prime\prime\prime}_{\mathrm{NJ}}(\lambda,X)<\frac{1}{2}\lambda \epsilon^{2}-(\epsilon-4)\lambda+1.
	\end{aligned}$$
	a contradiction. This completes the proof.
\end{proof}
Let  $\lambda \in(0,1)$  and  $a \in[0,1]$ . We define
$$
\begin{aligned}
	\phi(\lambda, a)=&\frac{1}{2} \sup \bigg\{\frac{1}{2}(\|\lambda x+(1-\lambda) a y\|^{2}+  \|\lambda x-(1-\lambda) a y\|^{2} \\
	& +\lambda(1-\lambda)\|x-a y\|^{2}+\lambda(1-\lambda)\|x+a y\|^{2}):x, y \in S_{X}\bigg\}.
\end{aligned}
$$

\begin{proposition} Let  $\lambda \in(0,1)$. The function  $\phi(\lambda, a) $ of the variable  $a$  is continuous on  $[0,1)$  and is convex and nondecreasing on  $[0,1]$.
\end{proposition}
\begin{proof}
Obviously, by the constant $C^{\prime\prime\prime}_{\mathrm{NJ}}(\lambda,X)$, we can easily get the result,
so we omit the proof.
\end{proof}
\begin{corollary}
Let  $X$  be a Banach space and  $\lambda \in(0,\frac{1}{2}]$ . Then
$$
C^{\prime\prime\prime}_{\mathrm{NJ}}(\lambda,X)= 2\phi(\lambda, a).
$$
\end{corollary}
The problem of this constant is that we don't know if there is a Banach space  $X$  with  $C^{\prime\prime\prime}_{\mathrm{NJ}}(\lambda,X)=1 $. However in this case we were able to calculate  $C^{\prime\prime\prime}_{\mathrm{NJ}}(\lambda,\ell_{q})$  for  $q \geq 2$ .
\begin{example}
Let  $q \geq 2 $ and  $\lambda \in(0,\frac{1}{2}]$ . Then
$$
C^{\prime\prime\prime}_{\mathrm{NJ}}(\lambda,\ell_{q})=\bigg(\frac{(\sqrt{1-\lambda}+\sqrt{\lambda})^{q}+(\sqrt{1-\lambda}-\sqrt{\lambda})^{q}}{2}\bigg)^{\frac{2}{q}}.
$$
\end{example}
\begin{proof}
	We will use Clarkson's inequality (see \cite{JAC}): for  $q \geq 2 $ and  $x, y \in \ell_{q}$ 
$$
	\|x+y\|^{q}+\|x-y\|^{q} \leq(\|x\|+\|y\|)^{q}+|(\|x\|-\|y\|)|^{q}.
$$
	
	Using this inequality and the definition of  $\phi(\lambda, a)$  we have that for $ \lambda \in(0,\frac{1}{2}], a \in   [0,1]$  and $ x, y \in \ell_{q}$  with  $\|x\|=\|y\|=1 $ :
	$$
	\begin{aligned}
		&\frac{1}{4}(\|\lambda x+(1-\lambda) ay \|^2+\|\lambda x-(1-\lambda) ay\|^{2}+\lambda(1-\lambda)(\|x-ay\|^{2}+\|x+ay\|^{2})) \\
		 \leq& \frac{1}{4}(2^{1-\frac{2}{q}}(\lambda+(1-\lambda) a)^{q}+|\lambda-(1-\lambda) a|^{q}+2^{1-\frac{2}{q}}\lambda(1-\lambda)((1-a)^{q}+(1+a)^{q}))\\
		=&\frac{2^{1-\frac{2}{q}}}{4}((\lambda+(1-\lambda) a)^{q}+|\lambda-(1-\lambda) a|^{q}+\lambda(1-\lambda)((1-a)^{q}+(1+a)^{q}))=g(a).
	\end{aligned}
	$$
	Hence,  $\phi(\lambda, a) \leq g(a) $. If one take  $x=(\frac{1}{2^{\frac{1}{q}}} , \frac{1}{2^{\frac{1}{q}}}, 0, \ldots)$  and  $y=(\frac{1}{2^{\frac{1}{q}}},-\frac{1}{2^{\frac{1}{q}}}. ,  0, \ldots  )$ we have the equality, thus  $\phi(\lambda, a)=g(a) $ and
	$$
	C^{\prime\prime\prime}_{\mathrm{NJ}}(\lambda,\ell_{q})=\sup _{a \in[0,1]} 2\phi(\lambda, a)=\sup _{a \in[0,1]} 2g(a).
	$$
	
	The maximum is attained at $ a_{0}=\sqrt{\frac{\lambda}{1-\lambda}} $, and  $g(a_{0})=\frac{1}{2}\bigg(\frac{(\sqrt{1-\lambda}+\sqrt{\lambda})^{q}+(\sqrt{1-\lambda}-\sqrt{\lambda})^{q}}{2}\bigg)^{\frac{2}{q}} $.
\end{proof}

\vspace{2em}

\noindent\textbf{Acknowledgments.}  Thanks to all the members of the Functional Analysis Research team of the School
of Mathematics and Physics of Anqing Normal University for their discussion and
correction of the difficulties and errors encountered in this paper. 

\noindent\textbf{Funding.} This research work was funded by Anhui Province Higher Education Science Research Project (Natural Science), 2023AH050487.

\noindent\textbf{Conflict of interest.} The authors declare that there is no conflict of interest in publishing the article.



\end{document}